\documentstyle[11pt,epsf]{amsart}

\newlength{\figboxwidth}
\newcommand{\bld}[1]{\medskip \noindent {\bf #1 }\nopagebreak}

\newtheorem{theorem}{Theorem}[section]
\newtheorem{prop}[theorem]{Proposition}

\newtheorem{cor}[theorem]{Corollary}

\newtheorem{conj}[theorem]{Conjecture}

\renewcommand{\qed}{\hspace*{\fill} 
\setlength{\unitlength}{1mm}
\begin{picture}(2.5,2.5)
      \put(0,0){\framebox(2.5,2.5){}}
  \end{picture}
\setlength{\unitlength}{1pt}}

\newcommand{\Id}{{\rm Id}}

\newcommand{\integers}{{\bf Z}}

\newcommand{\ratls}{{\bf Q}}
\newcommand{\reals}{{\bf R}}

\newcommand{\zed}{\integers}

\newcommand{\diam}{\operatorname{diam}}

\newcommand{\rank}{\operatorname{rank}}
\newcommand{\dist}{\operatorname{dist}}

\newcommand{\tauG}{\tau(G)}
\newcommand{\ub}{\underline{\rm{b}}}
\newcommand{\ue}{\underline{\rm{e}}}

\newcommand{\uw}{\underline{\rm{w}}}

\newcommand{\uom}{\underline{\omega}}

\newcommand{\tr}{\operatorname{tr}}

\newcommand{\Lap}{{\bold L}}

\newcommand{\cE}{{\cal{E}}}
\newcommand{\cC}{{\cal{C}}}
\newcommand{\cG}{{\cal{G}}}

\newcommand{\cM}{{\cal{M}}}
\newcommand{\cN}{{\cal{N}}}

\newcommand{\cS}{{\cal{S}}}
\newcommand{\cT}{{\cal{T}}}
\newcommand{\cX}{{\cal{X}}}

\newcommand{\tP}{{\tilde{P}}}

\title{On some extremal problems in graph theory}
\author{Dmitry Jakobson and Igor Rivin}
\thanks{The first author was partially supported by the 
NSF postdoctoral fellowship} 
\address{First author: 
IAS, School of Math, Princeton, NJ 08540 and 
Dept. of Mathematics 253-37, Caltech, Pasadena, CA 91125.  
Second author: Department of Mathematics, Warwick University\\
Coventry CV4 7AL, England 
and 
Dept. of Mathematics 253-37, Caltech, Pasadena, CA
91125.
} 
\date{\today}
\keywords{Graph, weighted, Laplacian, eigenvalue, cage, spanning tree, girth, 
diameter, mean distance}

\begin{document}
\bibliographystyle{jakobson}

\begin{abstract}
In this paper we are concerned with various graph invariants (girth, 
diameter, expansion constants, eigenvalues of the Laplacian, tree number) 
and their analogs for weighted graphs -- weighing the graph changes
a combinatorial problem to one in analysis.  We study both weighted and 
unweighted graphs which are extremal for these invariants.  In the 
unweighted case we concentrate on finding extrema among all (usually)
regular  graphs with the same number of vertices; we also study the
relationships  between such graphs.   
\end{abstract}
\maketitle

\section{Introduction}
\label{sec:intro}
In this paper we are concerned with various invariants of graphs related to 
distance (girth, diameter),  
expansion properties 
of the graph (Cheeger constant) and spectrum of the Laplacian 
(the first eigenvalue, the tree number).

The behavior of many of these invariants for random graphs has been 
extensively studied (\cite{Bollobas}).  For regular graphs, the results 
about weak convergence of graph spectra to the spectrum of the infinite 
regular tree (\cite{Kesten,McKay:81}) allow to get many 
results on the asymptotic behavior of spectral invariants (as the number 
of vertices increases).  

Often when certain bounds are established for the rate of growth of some 
invariant, one is interested in studying families graphs with
(asymptotically)  the {\em best possible} behavior with respect to a
given invariant (two  
well-known examples are Ramanujan graphs and cages).  

A related problem is to study the graphs which maximize (or minimize) 
a given invariant among all simple connected graphs with $n$ vertices 
and $m$ edges
(we denote the set of such graphs by $\cG_{n,m}$), or among $k$-regular 
connected simple graphs with $n$ vertices (we denote the set of such 
graphs by $\cX_{n,k}$).  
The sets $\cG_{n,m}$ and $\cX_{n,k}$ are finite, so maximizing 
a given invariant is a discrete optimization problem.  

On the other hand, $\cG_{n,m}$ (respectively, $\cX_{n,k}$) is naturally 
embedded in the complex $\cG_{n,m}^w$ of all {\em weighted} graphs with 
$n$ vertices and $m$ edges (respectively, the complex $\cX_{n,k}^w$ of
weighted $k$-regular graphs with $n$ vertices); the edge weights are 
normalized so that their sum is constant.  One can study the extremal 
problems for weighted analogs of graph invariants (see e.g. \cite{Fiedler}).
A natural question then is which {\em unweighted} simple graphs in 
$\cG_{n,m}$ or $\cX_{n,k}$ (all weights are 
equal to $1$) are local extrema in $\cG_{n,m}^w$ or $\cX_{n,k}^w$ 
for a given invariant.  In section \ref{sec:edges} we study this 
question for various graph invariants.  

In section \ref{sec:heuristics} we discuss some heuristics about the 
structure of graphs which have many spanning trees.  
In section \ref{sec:genustree} we give an upper bound on the rate of growth 
of the number of spanning trees for $6$-regular simple graphs on a fixed 
surface.  In section \ref{sec:lattice} we prove a result about 
the asymptotic properties of certain lattices, associated to 
graphs by Bacher, de la Harpe and Nagnibeda (cf. \cite{BHN}), 
for sequences of regular graphs.  Finally, in section 
\ref{sec:construct} we explain how the first $8$ cubic cages 
appear naturally in a sequence of graphs produced from $K_4$ by 
``greedy'' edge insertion with respect to the number of spanning 
trees.  

In the remainder of this section we review the definitions and some 
results about the asymptotic behavior of various graph invariants.  

\subsection{Definitions and notation}
The {\em adjacency} matrix $A$ of a graph $G\in \cG_{m,n}$ 
is a square matrix of size $n$  where $A_{i,j}$ is the number of edges 
joining the vertices $v_i$ and $v_j$.  We always consider loopless 
graphs, so $A_{ii}=0$; for simple graphs $A_{ij}$ 
is equal to $0$ or $1$.      

The nearest neighbor {\em Laplacian} $\Lap$ acts on functions on 
the set $V(G)$ of the vertices of $G$: given $f:V\to\reals$,  
$\Lap (f)(v)\ =\ \sum_{w\sim v} (f(v)-f(w))$.  Let $\Lap =\{L_{ij}\}$ be the 
matrix of $\Lap(G)$ of a (not necessarily simple) graph $G$; then $L_{ii}$ is 
the degree of the vertex $v_i$, and $-L_{ij}$ is the number of the edges 
joining $v_i$ and $v_j\neq v_i$ (equal to $0$ or $1$ for simple graphs).   
For $k$-regular graphs $\Lap=k\cdot\Id -A$.  

Let $A$ be the adjacency matrix of a graph $G\in \cX_{n,k}$, and 
let its spectrum (in the decreasing order) be given by
\begin{equation}\label{adj:spectrum}
   k=\lambda_1 >\lambda_2\geq\lambda_3\geq\ldots\geq\lambda_n\geq (-k)
\end{equation}
(the constant function is the eigenfunction with the eigenvalue $k$).  
The last eigenvalue $\lambda_n=-k$ if and only if $G$ is bipartite.  
The spectrum of $\Lap(G)$ is $\mu_j=k-\lambda_{j+1}, j=0,1,\ldots ,n-1$; 
it satisfies $0=\mu_0<\mu_1\leq\ldots\leq\mu_{n-1}\leq \min\{n,2k\}$.  

\subsection{Spectral properties of random graphs}
Alon and Boppana (\cite{Alon:86,LPS}) proved that 
\begin{equation}\label{liminf:adj}
   \lim_{n\to\infty}\; \left(\inf_{G\in \cX_{n,k}} \lambda_2(G)\right)
   \ \geq\ 2\;\sqrt{k-1}
\end{equation}
(cf. \cite{Nilli,Sole} for more precise results);   
$k$-regular graphs all of whose $\lambda_j\neq\pm k$ lie inside 
$[-2\sqrt{k-1},2\sqrt{k-1}]$ are called {\em Ramanujan} graphs 
(because of the relation to Ramanujan conjectures, cf. 
\cite{Sarnak:book,Lub:book}).  The 
first such graphs were constructed in \cite{LPS,Mar:88}; many
other constructions followed.   

\subsection{Complexity of a graph}
An important invariant of $G$ is the number $\tauG$ of spanning trees 
of $G$; it is sometimes called the {\em complexity} of $G$.
By Kirckhoff's theorem (\cite{Kirchhoff}), 
$$
n\; \tauG\ =\ \mu_1\mu_2\cdot\ldots\cdot\mu_{n-1}
$$
and $\tauG$ is equal to 
the determinant of any cofactor of the matrix $\Lap$ and 
(it is some times called the {\em determinant of Laplacian}).  
McKay proved in \cite{McKay:83} that 
\begin{equation}\label{limsup:trees}
   \lim_{n\to\infty}\; \left(\sup_{G\in \cX_{n,k}} (\tauG)^{1/n}\right)
   \ =\ \sigma_k\ =\ \frac{(k-1)^{k-1}}{(k(k-2))^{k/2-1}}
\end{equation}
and that a typical graph $G$ in $\cX_{n,k}$ satisfies 
$$
   \left(\tauG\right)^{1/n}\ \sim\ \sigma_k
$$
The above formula also holds for some Ramanujan graphs (\cite{Sarnak:ppt}). 

Alon in \cite{Alon:90} proved some results concerning 
\begin{equation}\label{liminf:trees}
   \delta_k\ =\ 
   \lim_{n\to\infty}\; \left(\inf_{G\in \cX_{n,k}} (\tauG)^{1/n}\right)
\end{equation}
Namely, he showed that for $k\geq 3$, $\delta_k\geq\sqrt{2}$ and that 
$\delta(k)/k=1+o(1)$ as $k\to\infty$.  Valdes in \cite{Valdes} determines 
the $2$-connected cubic graphs with the {\em minimum} number of spanning 
trees.  It follows from her results that if the limit in \eqref{liminf:trees}
is taken over all $2$-connected cubic graphs, the answer is $2^{3/4}$ 
(this was also proved by A. Kostochka in \cite{Kostochka}).  

\subsection{Extremals for girth}
We begin with girth.  It is easy to show that a $k$-regular graph 
($k\geq 3$) with girth $g\geq 3$ has at least $f_0(k,g)$ vertices where 
\begin{equation}\label{cage:low}
   f_0(k,g)\ =\ \left\{\aligned 
   \; &\ 1+k[(k-1)^{(g-1)/2} - 1]/(k-2),\\ 
   \; &\ 2[(k-1)^{g/2} - 1]/(k-2),\endaligned\right. 
   \qquad \left.\aligned
   \; & g\ odd\\
   \; & g\ even\endaligned\right.
\end{equation}
The $k$-regular graphs with girth $g$ and $f_0(k,g)$ vertices are called 
{\em Moore graphs} (sometimes they are also called {\em generalized polygons} 
for even $g\geq 4$); these graphs are unique if they exist.  
For $g=3$ those are complete graphs $K_{k+1}$, for 
$g=4$ those are complete bipartite graphs $K_{k,k}$.  It is proved in 
\cite[Ch. 23]{Biggs:book} that for $g>4$ Moore graphs only exist 
if $g=5$ and $k=3$ (Petersen graph), $k=7$ 
(Hoffmann-Singleton graph), and (possibly) $k=57$; or if $g=6,8$ or $12$.  
$k$-regular Moore graphs of girth $6$ exist if and only if a projective 
plane with $k$ points on the line does \cite{Singleton}.  

The $k$-regular graphs of girth $g$ which have the smallest possible number
of vertices are called $(k,g)$ {\em cages} (such graphs always exist, cf. 
\cite[Ch. 6]{HS:book}), and the number of their vertices is denoted by 
$f(k,g)$.  For $k\geq 3$ and $g$ such that the corresponding Moore graphs 
exist, those graphs are $(k,g)$ cages $f(k,g)=f_0(k,g)$; if Moore graphs 
don't exist, $f(k,g)>f_0(k,g)$.  
Very few cages other than Moore graphs are known (see \cite{Wong}, 
\cite[Ch. 6]{HS:book} and \cite[Ch. 23]{Biggs:book} for surveys of known 
results).  For cubic graphs, $f(3,g)$ is known for $3\leq g\leq 12$ 
(\cite{BMS}).  
For $k\geq 4$, the $f(k,7)$ is already unknown.  The best upper bounds for 
$f(k,g)$ come from various infinite series of symmetric graphs (cf.  
\cite{LPS,Mar:88,Mor,Chiu,BH,LUW}).   
It is known that 
\begin{equation}\label{girth:lim}
   1/2\ \leq\ 
\liminf_{g\to\infty} \frac{\log_{k-1} f(k,g)}{g}\ \leq\ 
\limsup_{g\to\infty} \frac{\log_{k-1} f(k,g)}{g} 
   \ \leq\ 3/4
\end{equation}
For all known infinite families of regular graphs, 
$$
\limsup_{g\to\infty}(\log_{k-1}f(k,g))/g=3/4
$$ 
it is sometimes conjectured 
that the upper bound in \eqref{girth:lim} is actually an asymptotic result.  

\subsection{Extremals for complexity}
We next turn to the number of spanning trees.  It is easy to see that the 
complete graph on $n$ vertices has the maximum number of trees among all 
multigraphs with the same number of vertices and edges.  
In a series of papers (see
e.g. \cite{CK,Kel:74,Kel:76,Kel:80,Kel:96,Kel:97}) Kelmans
characterized the  
graphs, obtained by removing a small number of edges form the complete 
graph, which have the maximum number of spanning trees among all the graphs 
with the same number of vertices and edges.  Cheng (\cite{Cheng}) proved 
that the complete multipartite regular graphs also have this property; 
see also \cite{Const}.  

We notice that $K_{k+1}$ and $K_{k,k}$ are, respectively, the unique $(k,3)$ 
and $(k,4)$ cages.  One can also check (see \cite{Valdes}) that the Petersen 
graph (which is the $(3,5)$ cage) has the most trees ($2000$) among 
all cubic simple graphs with $10$ vertices.  We conjecture
\begin{conj}\label{conj:Moore:trees}
Moore graphs have the maximum number of spanning trees among all simple 
graphs with the same number of vertices and edges.
\end{conj}
Some evidence supporting this conjecture is discussed in section 
\ref{sec:construct}.  

{\bf Acknowledgements.} The authors would like to thank N.~Alon,  
P.~Sarnak, P.~Seymour and R.~Wilson for helpful conversations, and A.~Kelmans 
for explaining many 
results in extremal graph theory during his visit to Caltech in the spring 
of 1997, as well as for suggesting many improvements to the manuscript.  
The first author would like to thank the Institute for Advanced Study at 
Princeton for their hospitality in the fall of 1997 and in the spring of 1998.


\section{Local extrema for weighing the edges}
\label{sec:edges}  

In this chapter we shall discuss several extremal problems for weighted 
graphs.  We will be concerned with finding extrema for various functionals 
on the {\em weighted edge polytope} $P$ (which is actually a simplex) 
defined below.  Let $G$ be a simple 
connected graph, and let $E(G)=\{e_1,\ldots ,e_m\}$ be the set of its 
edges.  Then 
\begin{equation}\label{weight:sum}
P(G)\ :=\ \left\{(x_1,\ldots ,x_m):\; x_j\geq 0,\;\sum_{j=1}^m x_j\ =
\ m\right\}
\end{equation}
where the (nonnegative) real number $x_j$ is called the {\em weight} 
of $e_j$; the corresponding graph is denoted by $G(x_1,\ldots,x_m)$.    
The graph $G(1,1,\ldots,1)$ is called an {\em unweighted graph}.  
We introduce the following notation:
$$
\uom\ =\ (1,\ldots,1).
$$  

The adjacency matrix of $G(x_1,\ldots,x_m)$ is defined by 
$$
A_{ij}\ =\ \begin{cases} 0, &  i=j\ or\ v_i\not\sim v_j \\
                     x_k, & e_k=(v_i v_j) \end{cases}
$$
The {\em degree} $\deg(v_i)$ of a vertex $v_i$ is defined by 
$$
\deg(v_i)\ =\ \sum_{j:(v_iv_j)=e_k\in E(G)} x_k
$$
The Laplacian is defined by 
\begin{equation}\label{Lap:defn}
L_{ij}\ =\ \begin{cases} deg(v_i), &  i=j \\ 
                         0, & v_i\not\sim v_j \\
                         -x_k, & e_k=(v_i v_j) \end{cases}
\end{equation}
The notation for the eigenvalues of $A$ and $\Lap$ will be the same as 
for the unweighted graphs.  The proof of Kirchhoff's theorem extends to 
weighted graphs (\cite{Chung}); we shall continue to call the determinant 
$\tauG$ of any cofactor of $L_{ij}$ the (weighted) tree number of $G$.  

One can regard $\deg(v_i)$ as a combinatorial analog of a curvature at 
$v_i$.  Motivated by an analogy with metrics of constant curvature, we 
define another polytope $\tP(G)\subset P(G)$ by 
\begin{equation}\label{degree:preserve}
\tP(G)\ :=\ \left\{(x_1,\ldots ,x_m):\; x_j\geq 0,\;\deg(v_i)=d_i\; 
\forall i\right\}
\end{equation}
where $d_i=\#\{j\neq i:v_i\sim v_j\}$ is the degree of $v_i$ in the 
unweighted graph $G$.

We will consider the problem of maximizing or minimizing various 
functionals over $P$ and $\tP$.  Some of the functionals (girth and 
diameter) will be linear in $x_j$-s, and so we will be solving a 
problem in linear programming (\cite[Ch. 7]{Schrijver}).  Other 
functionals (expansion constants) will not be linear in $x_j$-s on $P$, 
but will become linear on $\tP$.  

We shall be often concerned with the following question:

{\bf Question:} When is 
$\uom=(1,1,\dots ,1)$ an extremum (in $P$ or $\tP$)?
   
In other words, when is the unweighted graph an extremum for a 
given functional?  In some cases we are able answer the above question
completely.

\subsection{Girth and diameter} 
We turn to ths study of the girth and the diameter.  
We make several definitions:
The {\em length} of a path in a weighted graph $G$ is the 
sum of the weights of the edges in the path.  The girth $\gamma(G)$ is the 
length of the shortest cycle in $G$.  The {\em distance} 
between two vertices of $G$ is the length of the shortest path connecting 
them; the {\em diameter} $\diam(G)$ is equal to the largest such distance.  

We are interested in the points $x\in P(G)$ where $\gamma(x)$ takes 
its maximal value and where $\diam(x)$ takes its minimal value on $P$. 
Denote by $C_1,\ldots,C_s$ all the cycles of $G$.  Then finding the 
maximum of the girth leads to solving the following problem in linear 
programming: 
$$
\max_{x\in P}\left(\min_{1\leq i\leq s}\sum_{e_j\in C_i} x_j\right).
$$
We assume that any edge on a cycle is visited at most once, so there are 
finitely many cycles; we make a similar assumption about the paths in 
the next definition.

Similarly, for two vertices $u,v$ of $G$ denote by 
$\{\pi_1(u,v),\ldots,\pi_\sigma(u,v)\}$  the set of all paths in $G$ 
connecting $u$ and $v$ (where $\sigma=\sigma(u,v)$).  Finding the 
minimum of the diameter leads to the following problem:
$$
\min_{x\in P}\left(\max_{u,v\in G}\min_{1\leq j\leq\sigma(u,v)}
\sum_{e_k\in \pi_j} x_k\right).
$$
 
We now fix some $x\in P(G)$ and consider the corresponding graph 
$G=G(x)$.  
Two vertices $u, v$ of $G$ are called {\em diametrically opposite} if 
$\dist(u,v)=\diam(G)$.  Any path connecting two diametrically opposite 
points which has length $\diam(G)$ is called a {\em meridian}.  
Also, any shortest cycle  is called a {\em systole}.

To study how girth changes, we remark that for small weight changes 
the new systoles will come from among the weighted old systoles.
Denote the systoles of $G$ by $S_1,S_2,\ldots ,S_r$, 
and form the ($m$-by-$r$) edge--cycle incidence matrix $S_{ij}$.  
Turning to the diameter, we remark that for small weight changes the 
new meridians will come from among the weighted old meridians.  
Denote the meridians of $G$ by $M_1,M_2,\ldots ,M_t$ and form the 
($m$-by-$t$) edge-meridian incidence matrix $M_{ij}$.  
Denote by $\cS$ (respectively, $\cM$) the cone formed by all linear 
combinations of column vectors of $S_{ij}$ (respectively, $M_{ij}$) 
with nonnegative coefficients.  Clearly, both $\cC$ and $\cM$ lie 
inside the first quadrant.  

\begin{theorem}
\label{thm:girdiam}
   The graph $G$ is a local maximum for girth (respectively, a local 
   minimum for diameter) if and only if $\uom$ lies in $\cS$ 
   (respectively, $\cM$).
\end{theorem}

\bld{Proof.}
Let $G(x^0)$ be the graph in question, and let $x^1$ be another point in 
$P(G)$.  Denote by $y$ the vector $x^1-x^0$ connecting the points in $P$. 
The graph $G(x^0)$ is a local maximum for the girth if and only if the 
inequality 
$$
S\; y\ >\ 0
$$
has no solutions for small enough $y$.  
By Farkas' lemma (\cite[Ch. 7]{Schrijver}), this condition is  
equivalent to the condition $\uom\in\cS$.  

Similarly, the graph $G(x^0)$ is a local maximum for the girth if and 
only if the inequality 
$$
M\; y\ <\ 0
$$
has no solutions for small enough $y$.  
Again by Farkas' lemma, this is equivalent to the condition $\uom\in\cM$.  

\qed

An obvious sufficient condition for the assumptions of the Theorem 
\ref{thm:girdiam} to hold is for every edge to be contained in the
same number  of systoles (diameters), so that edge-transitive graphs are 
always satisfy this assumption.  Another large class of graphs where 
every edge is contained in the same number of systoles consists of graphs 
arising from $4$-vertex-connected proper triangulations of surfaces
(every edge is contained in exactly two triangles). 

On the other hand, the assumptions of Theorem \ref{thm:girdiam} are clearly 
{\em not} satisfied when there is an edge in $G$ which is not contained in 
any systole (diameter).

\subsection{Expansion constants}

We now turn to Cheeger and expansion constants.  Let $G$ be an unweighted 
graph.  Given a partition of its vertex set $G=A\cup B$ into {\em disjoint} 
sets $A$ and $B$, define $E(A,B)$ 
to be the set of all edges with one endpoint in $A$ and another in $B$.  
Let $|E(A,B)|$ be the sum of the  weights of the edges in 
$E(A,B)$,\footnote{It is equal to the number of the edges in $E(A,B)$ for 
unweighted graphs.} 
and let $|A|$ denote 
the sum of degrees of all vertices in $A$ (also called the {\em volume} of 
$A$).\footnote{Often $|A|$ denotes the number of vertices in $A$, which for 
regular graphs is proportional to our definition; our definition allows a 
straightforward extension to weighted graphs.}  
The expansion constant $c(G)$ is defined by 
\begin{equation}\label{def:exp:const}
   c(G)\ =\ \min_{G=A\cup B}\; \frac{|E(A,B)|}{|A|\cdot |B|}
\end{equation}
The {\em Cheeger constant} $h(G)$ is defined by 
\begin{equation}\label{def:Ch:const}
   h(G)\ =\ \min_{G=A\cup B}\; \frac{|E(A,B)|}{\min(|A|,|B|)}
\end{equation}

To prove an analog of Theorem \ref{thm:girdiam} for the above constants, 
we find it convenient to maximize Cheeger and expansion constants over 
$\tP(G)$ instead of $P(G)$.  The reason is that for $A\subset V(G)$ the 
volume $|A|$ does not depend on $x\in\tP(G)$ (since the degrees are 
preserved), and so $c(G)$ and $h(G)$ become linear functionals on $\tP(G)$. 

Let $N_{ij}$ denote the $m$-by-$n$ edge--vertex incidence matrix.  
The dimension of $\tP$ is equal to $m-\rank(N_{ij})$.  
We denote by $\cN$ the cone spanned by linear combination of the column 
vectors of $(N_{ij})$ with nonnegative coefficients.\footnote{Here $\cN$ 
stands for {\em normalization} conditions; in general, the more conditions we 
impose, the larger the cone $\cN$ is; vector $\uom$ should always be 
in $\cN$.}  

To state the next theorem, we need to define several matrices.  Given 
$x\in\tP$, let $G(x)$ be the corresponding graph.  
An edge cut $H$ of $G(x)$ is called an {\em expansion edge 
cut} (respectively, a {\em Cheeger edge cut}) if it is one of the cuts for 
which the expansion constant (respectively, the Cheeger constant) of $G$ 
is attained.  Let $H_1,H_2,\ldots ,H_r$ (respectively, $H^1,H^2,\ldots ,H^t$)
be the expansion (respectively, Cheeger) edge cuts.  
 
Denote by $A(H)$ and $B(H)$ the disjoint subsets into which $H$ separates 
$V(G)$.  We now define two matrices: an $m$-by-$r$ {\em expansion} matrix 
$(E_{ij})$ and an $m$-by-$t$ {\em Cheeger matrix} $(C_{ij})$ by 
\begin{equation}\label{matrix:exp}
   E_{ij}\ =\ \left\{\aligned 
   \ & 1/(|A(H_j)|\cdot |B(H_j)|), \qquad e_i\in H_j\\
   \ &\ 0, \qquad otherwise.\endaligned\right.
\end{equation}
and by 
\begin{equation}\label{matrix:Cheeger}
   C_{ij}\ =\ \left\{\aligned 
   \ & 1/\min(|A(H_j)|,|B(H_j)|), \qquad e_i\in H^j\\
   \ &\ 0, \qquad otherwise.\endaligned\right.
\end{equation}
Let $\cE$ (respectively, $\cC$) denote the cone formed by linear 
combinations of column vectors of $E_{ij}$ (respectively, $C_{ij}$).  

We are now ready to formulate and prove
\begin{theorem}
\label{thm:Chexpan}
   A graph $G$ is a local maximum for the expansion constant 
   (respectively, the Cheeger constant) for degree-preserving 
   edge valuations if and only if $\cE\cap\cN$ (respectively, 
   $\cC\cap\cN$) contains a nonzero vector. 
\end{theorem}

\bld{Proof}: The proof is similar to that of Theorem
\ref{thm:girdiam}. 
\qed

We remark that edge-transitive graphs satisfy both 
assumptions of Theorem \ref{thm:Chexpan} (for such graphs 
$\uom\in(\cE\cap\cN)$ and $\uom\in(\cC\cap\cN)$).  
Also, as in the case of girth and diameter the assumptions 
are not satisfied if some edge is not contained in any 
Cheeger (respectively, expansion) edge cut.

\subsection{The tree number}
By the weighted version of Kirckhoff's theorem (\cite{Chung})   
\begin{equation}\label{det:tree-exp}
    \tauG\ =\ \sum_{T\in\cT(G)} \; \prod_{e_j\in T} x_j
\end{equation}
where the sum is taken over the set $\cT(G)$ of the spanning trees of $G$.  

Finding the maximum of $\tauG$ on $P$ becomes a Lagrange multiplier 
problem with the constraints given by Eq. (\ref{weight:sum}). 
The condition for $x\in P(G)$ to be a critical point for $\tauG$ is 
\begin{equation}
\frac{\partial\tauG}{\partial x_1}\ =\ 
\frac{\partial\tauG}{\partial x_2}\ =\ \ldots\ =\ 
\frac{\partial\tauG}{\partial x_m}
\end{equation}

The partial derivatives above are given by 
$$
\tau_j\ =\ \frac{\partial\tauG}{\partial x_j}\ =\ \sum_{T\in\cT(G)}
\prod_{k\neq j;e_k\in T} x_k.  
$$ 
The ratio $\tau_j/\tauG$ is called the {\em effective resistance} of $e_j$. 

We have thus proved:
\begin{prop}\label{prop:detcrit}
The graph $G(x)$ is a critical point for $\tauG$ if and only if the effective 
resistances of all edges are the same.
\end{prop}

If an unweighted graph satisfies the assumptions of 
Proposition \ref{prop:detcrit} then every edge of this graph 
is contained in {\em the same} number of spanning trees.  Such 
graphs were studied by Godsil in \cite{Godsil}; he calls these 
graphs {\em equiarboreal}.  Obviously, all edge-transitive graphs 
(the automorphism group acts transitively on the edges) are 
equiarboreal.\footnote{See \cite{Bouwer} for examples of edge-transitive 
graphs which are not vertex-transitive.} 
Godsil gives several more sufficient conditions for 
a graph to be equiarboreal; in particular, any distance-regular graph 
and any color class in an association scheme is equiarboreal (the least 
restrictive condition Godsil gives is for a graph to be 
{\em 1-homogeneous}).  By an easy counting argument one can show that 
for an unweighted equiarboreal graph 
\begin{equation}\label{equiarb:effres}
   T_1\ =\ T_2\ =\ \ldots\ =\ \tauG\cdot (n-1)/m
\end{equation}
(this is actually the result of Foster, cf. \cite{Foster})
so the necessary condition for a graph to be equiarboreal is 
$$
   m\ |\ (n-1)\tauG.
$$

We remark that the graphs which have the most spanning trees among the 
regular graphs with the same number of vertices are not necessarily 
equiarboreal, and vice versa.  For example, the $8$-vertex M\"obius wheel 
(cf. section \ref{sec:construct}) which has the most spanning trees among 
the $8$-vertex cubic graphs is not equiarboreal (cf. also \cite{Valdes}),   
while the {\em cube} (which is certainly edge-transitive, hence 
equiarboreal) has the {\em second biggest} number of spanning trees 
among the $8$-vertex cubic graphs.  

We next study the function $\tauG$ on $P(G)$.  Surprisingly, it turns
out to be concave, so there is at most one  critical point, and if
there is one it has to be a global maximum. The eigenspace $E_0$
corresponding to the eigenvalue $\mu_0=0$ of $\Lap$  
on $G(x)$ is the same for all $x$ in the interior of $P$ (and consists 
of constant functions). Let $\Pi_0$ be the operator of projection onto 
$E_0$; its matrix is given by $(1/n)\; J$, where $J$ is the matrix of 
$1$-s.  The eigenvalues of 
$$
M\ =\ \Lap\; +\; \Pi_0
$$ 
are $\{1,\mu_1,\ldots ,\mu_n\}$,   
hence   
\begin{equation}\label{Lap:posdef}
n^2\;\tauG\ =\ \det(J+\Lap).
\end{equation} 
This result is due to Temperley (\cite{Temperley}).
Thus it suffices to study the concavity of the positive-definite 
operator $M$.  

\begin{theorem}\label{logdet:concave}
The function $\tauG$ is concave on (the interior of) $P$.
\end{theorem}

The theorem immediately implies 
\begin{cor}
If $z$ is a critical point of $\tauG$ in the interior of $P$, then it
is the global maximum for $\tauG$.  In particular, the unweighted
graph is the global maximum for equiarboreal graphs.
\end{cor}

\bld{Proof of Theorem \ref{logdet:concave}.}
It suffices to prove that the determinant is a concave function on the set
of positive definite symmetric matrices.  Let $Q$ be such a matrix, and 
let 
$$
Q(t)\ =\ Q+tB,\ \ t\in\reals
$$
be a line of symmetric matrices through $Q$.  Then 
$$
\frac{d\log\det (Q(t))}{dt}\ =\ \tr(BQ^{-1}), 
$$
and 
$$
\frac{d^2\log\det (Q(t))}{dt^2}\ =\ -\tr(BQ^{-1}BQ^{-1}). 
$$

It suffices to show that the last trace is strictly positive.  
The matrix $R=Q^{-1}$ is positive definite.   Conjugate $R$ to a  
diagonal matrix with eigenvalues $\uw=(w_1,\ldots,w_n)$.  Then
$$
\tr(BRBR)\ =\ \uw S\; \uw^t
$$
where $S_{ij}=R_{ij}R_{ij}$ is the {\em Hadamard product} of $R$ with 
itself.  By Schur product theorem (\cite[Theorem 7.5.3]{HJ:book}), $S$ 
is positive definite, which finishes the proof.  
\qed

We remark that Theorem \ref{logdet:concave} can be considered as an 
analogue of the convexity of the determinant of Laplacian in a 
conformal class of metrics on surfaces of genus greater than or 
eqaul to one, proved by Osgood, Phillips and Sarnak in \cite{OPS}.  
It would be interesting to characterize graphs $G$ where $\tauG$ achieves 
its maximum in the interior of $P(G)$.  

\subsection{Eigenvalues of the Laplacian}

We introduce some notation which will be used in studying the 
Taylor expansion of the eigenvalues. 

Given two points $x^0,x^1\in P(G)$, let $y=x^1-x^0$ and let 
\begin{equation}\label{newts}
x(t)\ =\ x^0\; +\; ty,\ \ 0\leq t\leq 1
\end{equation}
be the segment (call it $\pi(t)$) connecting the two points in $P$.  
The normalization condition reads 
\begin{equation}\label{norm:wts}
   y\cdot\uom^t\ =\ 0
\end{equation} 
The basis of all $y$-s satisfying \eqref{norm:wts} is given by vectors 
\begin{equation}\label{basewts}
   \ub_{ij}\ =\ \ue_i-\ue_j
\end{equation}
for $i\neq j$ where $\ue_i$-s are the standard basis vectors.  

We next turn to studying how the eigenvalues of the Laplacian 
depend on $t$.  We look at one parameter family $\Lap(t)$ 
of Laplacians on $G$ corresponding to the points on $\pi(t)$.  
The Laplacian is given by the matrix 
$$
   \tilde{\Lap}(t)\ =\ \Lap\; +\; tB
$$ 
By standard results in perturbation theory (\cite{Kato}), there 
exists an orthonormal basis of the eigenvectors of $Q$ such that the 
the eigenvalues and the eigenvectors admit Taylor expansions in $t$.

Let $\mu=\mu(0)>0$ be an eigenvalue of the multiplicity $p\geq 1$. 
Let $E_\mu$ be the corresponding ($p$-dimensional) eigenspace, and 
let $P_\mu$ denote the corresponding projection.  The 
condition for the graph $G$ to be a critical point for $\mu$
is that the matrix $P_\mu BP_\mu$ be indefinite for all $B$.  
We now compute the eigenvalues of $P_\mu BP_\mu$ for the edge 
weightings given by \eqref{newts} with $\ub$ as in \eqref{basewts}.  
We distinguish between two cases: 
\begin{itemize} 
\item[a)] edges $e_1$ and $e_2$ are disjoint.
\item[b)] edges $e_1$ and $e_2$ share an endpoint
\end{itemize}

We first consider the case a). Let $f_1,f_2,\ldots ,f_p:V\to\reals$ be a  
basis of eigenvectors of $E_\mu$.  For convenience, we denote 
$f_j(v_1)-f_j(v_2)$ by $z_j$ and $f_j(v_3)-f_j(v_4)$ by $w_j$.  
The eigenvalues of $P_\mu BP_\mu$ are given by the eigenvalues of the 
linear form $F(\alpha_1,\ldots ,\alpha_p)$ defined by 
\begin{equation}\label{eigen:prod}
F(\alpha_1,\ldots ,\alpha_p)\ =\ \langle 
B(\alpha_1f_1+\ldots +\alpha_pf_p), f_j\rangle
\end{equation}
It is easy to see that the matrix of $F=F_{ij}$ is given by 
$$
   F_{ij}\ =\ z_iz_j-w_iw_j
$$
The eigenvalue $\mu$ is critical if $F$ is indefinite for any choice 
of $B$.  

We remark that $F$ is a difference of two rank one matrices, 
$F=z\otimes z -w\otimes w$ where $z=(z_1,\ldots ,z_p)$ and 
$w=(w_1,\ldots ,w_p)$ and so has rank at most two.  Computing the 
coefficient $C_{p-2}$ of $t^{p-2}$ in the characteristic polynomial of $F$,  
we find that 
$$
C_{p-2}\ =\ -\sum_{1\leq i<j\leq p} (z_iw_j-z_iw_j)^2.
$$
Therefore, $C_{p-2}\leq 0$ and is strictly negative unless $x$ and 
$y$ are proportional or one of them is zero.  Now, if $C_{p-2}<0$ 
then $F$ is indefinite (since it has at most two nonzero eigenvalues).  
If $C_{p-2}=0$ (and $F$ has rank one) and if $\mu$ is critical then 
$F\equiv 0$.  

Consider first the case of a simple eigenvalue ($p=1$).  Then $\mu$ is 
critical if $F=F_{11}=0$ for every choice of $B$ and therefore 
$z_1^2-w_1^2\ =\ 0$ or $z_1=\pm w_1$.   Recalling the definition of 
$z_1$ and $w_1$, we see that for the eigenvector $f_1=f$
\begin{equation}\label{df:equal}
|f(v_1)-f(v_2)|\ =\ |f(v_3)-f(v_4)|
\end{equation} 
The following proposition is proved by repeating the argument for all 
pairs of nonadjacent edges; if some edge $e$ of $G$ is adjacent to every  
other edge, we prove an analogue of \eqref{df:equal} for the case b) 
of adjacent edges (it was proved before for the case a) of 
non-adjacent edges).   
\begin{prop}\label{eigen:crit}
A simple eigenvalue $\mu$ of the Laplacian on $G$ is critical
if and only if the corresponding eigenvector $f:V\to\reals$ satisfies 
\begin{equation}\label{df:const}
   |f(u)-f(v)|\ \equiv\ const
\end{equation}
\end{prop}

We now study the graphs which admit an eigenvector $f:V\to\reals$ 
satisfying \eqref{df:const} for some $c\geq 0$.  
If $c=0$ then $f$ is a multiple of a constant vector and 
so has eigenvalue zero which is a contradiction. 
If $c\neq 0$ then it is easy to see that the graph $G$ cannot 
have odd cycles and hence is bipartite. Namely, let 
$u_1u_2\ldots u_l$ be a cycle.  Then (putting $u_l=u_0$) 
$\sum_{i=1}^l (f(u_i)-f(u_{i-1}))=0$.  But each term in the 
sum is equal to $\pm c$, and since the number of terms in the sum is
odd, they cannot add up to $0.$

We now want to study the unweighted $k$-regular graphs which have an 
eigenvector (corresponding to an eigenvalue $\mu>0$) satisfying 
\eqref{df:const} (without necessarily assuming that $\mu$ is simple).    
We shall rescale the eigenvector so that $c=1$ in 
\eqref{df:const}.  From \eqref{Lap:defn} and \eqref{df:const} it follows 
that for each vertex $u$ the expression $\mu\cdot f(u)$ can only take 
one of the values $k, k-2,k-4,\ldots, -k+2,-k$.  
Consider first the vertex $u_0$ where $f(u)$ takes its maximal value $a$
(by changing the sign if necessary we can assume that $a>0$).  
It follows that $f$ takes value $a-1$ on all the neighbors of $u_0$, 
hence 
$$
   a\mu\ =\ k
$$ 

Next, consider any neighbor $u_1$ of $u$.  The value of $f$ at any 
neighbor of $u_1$ can be either $a$ (let there be $r_1\geq 1$ such 
neighbors; $u_0$ is one of them); or $a-2$ (it follows that there are 
$k-r_1$ such neighbors).  From \eqref{Lap:defn} it follows that 
$$
   \mu(a-1)\ =\ k-2r_1
$$
It follows from the last two formulas that 
\begin{equation}\label{eig:even}
   \mu\ =\ 2r_1
\end{equation}
where $r_1\geq 1$ is a positive integer.  If $r_1=k$, then $\mu=2k$ is 
the largest eigenvalue of $\Lap$.  

We next define the {\em level} of a vertex $u$ to be equal to $j$ if 
$f(u)=a-j$; we denote the set of all vertices of $G$ at level $j$ 
by $G_j$.  It is easy to see that if $u\in G_j$ has $r_j$ neighbors 
where $f$ takes value $a-j+1$ then
$$
   \mu(a-j)\ =\ k-2r_j
$$ 
It follows that $r_j$ is the same for all $u\in G_j$.
Using \eqref{eig:even} we see that 
$$
   r_1\cdot j\ =\ r_j
$$
Consider now a ``local minimum'' $u\in G_N$.   Then $r_N=k$, and we see 
that 
\begin{equation}\label{eig:divides}
   r_1\; |\; k
\end{equation} 

Let $n_j$ denote the number of vertices in $G_j$.  Counting the vertices 
connecting $G_j$and $G_{j+1}$ in two different ways, we see that 
for all $0\leq j\leq N-1$, 
$$
n_j(k-r_j)\ =\ n_{j+1}r_{j+1}
$$
Consider the case $r_1=1,\mu=2$.  It follows from the previous 
calculations that $r_j=j$ and that $N=k$.  Accordingly, 
$n_j=n_0{k\choose j}$ and 
\begin{equation}\label{vertn:dfconst:2}
|G|\ =\ 2^k\; n_0
\end{equation}
We next describe a class of graphs admitting an eigenvector of $\Lap$ with 
$\mu=2$ satisfying \eqref{df:const}.

An obvious example of such a $k$-regular graph is the $k$-cube, and any  
such graph has the same number of vertices as a disjoint union of $n_0$ 
cubes by \eqref{vertn:dfconst:2}.  Start now with such a union, choose 
the partition of the vertices of each cube into ``levels'' and take 
two edges $u_1u_2$ and $u_3u_4$ in two different cubes such that 
$u_1,u_3$ are both in level $j$ while $u_2,u_4$ are both in level $j+1$.
If we perform an edge switch 
$$
(u_1u_2),(u_3u_4)\ \to\ (u_1u_4),(u_3u_2)
$$ 
then the number of the connected components of our graph will decrease 
while the eigenvector $f$ will remain an eigenvector with the same 
eigenvalue.

Performing sequences of edge switches as described above, we obtain 
examples of connected graphs satisfying \eqref{df:const} and 
\eqref{vertn:dfconst:2} for any $n_0$.  Conversely, it is easy to show 
that starting from a graph satisfying \eqref{df:const} and 
\eqref{vertn:dfconst:2} and having chosen a partition of its vertices 
into levels one can obtain $n_0$ disjoint $k$-cubes by performing 
a sequence of edge switches as above.   

We now want to consider the case when $\mu=\mu_1$ is the lowest eigenvalue 
of the Laplacian.  The first remark is that then necessarily $\mu\leq k$, 
and $\mu=k$ only if $G=K_{k,k}$.  Next, we want to consider ``small'' 
$k$ for which $k-2\sqrt{k-1}$ (the ``Ramanujan bound'') is less than $2$
(this happens for $3\leq k\leq 6$).  
It then follows from the results of Alon (\cite{Nilli}) that the diameter 
of $G$ (and hence the number of vertices in $G$) is bounded above.  
\begin{prop}\label{mu1:finite}
For $3\leq k\leq 6$ there are finitely many $k$-regular graphs for 
which the condition \eqref{df:const} is satisfied for an eigenvector 
of $\mu_1$. 
\end{prop}   

We next discuss graphs which have an eigenvector satisfying \eqref{df:const} 
with the eigenvalue $\mu=2r_1>2$.  Recall that by \eqref{eig:divides} $r_1|k$. 
By counting the edges connecting the vertices in two consecutive levels one 
can show (as for $\mu=2$) that the number of vertices satisfies 
\begin{equation}\label{vertn:general}
|G|\ =\ 2^{(k/r_1)}\; n_0
\end{equation}
Also, since any vertex $u_1\in G_1$ has $r_1$ distinct neighbors in $G_0$, 
$$
n_0\geq r_1.
$$
It is easy to construct examples of regular graphs which have eigenvectors 
with the eigenvalue $\mu>2$ satisfying \eqref{df:const}; the construction 
is similar to that for $\mu=2$.  

We summarize the previous results: 
\begin{theorem}\label{eigencrit:character}
Let $G$ be a $k$-regular graph which has an eigenvector of $\Lap$ with an 
eigenvalue $\mu$ satisfying \eqref{df:const}.  Then $G$ is bipartite, 
$\mu=2l$ is an even integer dividing $2k$, the number of vertices of 
$G$ is divisible by $2^{(k/l)}$, and for $n_0\geq l$ there exist such 
graphs with $n=2^{(k/l)}n_0$ vertices.  
\end{theorem}


\section{Some heuristics for the number of spanning trees}
\label{sec:heuristics} 

Here we shall give some heuristic arguments concerning some properties 
that graphs with many spanning trees are expected to have.  
One such property is having large girth; the sequences of cubic graphs 
described in section \ref{sec:construct} seem to support that expectation. 
We shall give some heuristic arguments in its support.   

We start with the following formula from \cite[p. 153]{McKay:83} for the number 
$\tauG$ of spanning trees of a $k$-regular graph $G$ with $n$ vertices:  
\begin{equation}\label{tree:girth}
\tauG\ =\ kI_k(t)\lim_{t\to 0_+}\left(\frac{n}{n-1} + 
\frac{1}{(n-1)I_k(t)}\sum_{r=3}^\infty
{t\choose r}(-1)^r u_r k^{-r}\right)^{1/t} ,
\end{equation}
where $I_k(t)$ is a function depending on the valency $k$ only (and 
hence {\em the same} for all $k$-regular graphs), while $u_r$ is the 
number of walks which are not {\em totally reducible} (i.e. which don't 
``come'' from the infinite tree $T_k$).  
In particular, $u_r=0$ for $r<\gamma(G)$.  

The coefficients $(-1)^r{t\choose r}$ in the 
sum in \eqref{tree:girth} are always negative for small $t$, so each 
nonzero term in the sum {\em decreases} the value of the expression in 
brackets in \eqref{tree:girth}.  
Heuristically then, the more $u_r$-s are equal to zero, the larger 
the expression in brackets is in \eqref{tree:girth} and so the larger 
is the limit which is $\tauG$.   
Now, if we want to compare $\tauG$ for different graphs in $\cX_{n,k}$ 
those with larger girth $\gamma(G)$ have more $u_r$-s equal to zero and 
therefore are the natural candidates for large values of $\tauG$, and even 
more so Moore graphs, since they are 
the unique graphs in $\cX_{n,k}$ having the largest girth.  

Of course, we have only looked at the first 
nonzero term in the infinite series, but we remark that the subsequent 
$u_r$-s get multiplied by $(-1)^r{t\choose r}k^{-r}$ which is 
exponentially decreasing in $r$.  This seems to give at least a partial 
heuristic explanation of the cages appearing in section \ref{sec:construct}.  

It also seems reasonable to conjecture that graphs with higher connectivity 
have more spanning trees.  Indeed, one can show that if an edge $e$ 
is a bridge in a one-connected graph $G=G_1\cup\{e\}\cup G_2$ and if 
$G_1$ (say) is two-connected, then one can make an edge switch in $G$ 
which will increase the number of spanning trees (cf. \cite{Kel:97} for
the discussion of operations on graphs increasing their number of 
spanning trees).  However, it seems to be more difficult to 
prove that for a two-connected graph there exists an edge switch which 
simultaneously 
increases the number of spanning trees.  We believe this to be true in
essentially all cases, however.


\section{Graphs on a fixed surface}
\label{sec:genus} 

\subsection{The tree number for regular graphs on a surface}
\label{sec:genustree}
In this section we shall give an upper bound for the number of spanning trees 
for $6$-regular graphs on a fixed surface $S$ which is strictly smaller than 
the McKay's bound \eqref{limsup:trees} for general $6$-regular graphs.   

Let $G$ be a $k$-regular simple graph on $n$ vertices embedded in an 
orientable surface $S_g$ of genus $g$.  Let $V,E,F$ denote the number of 
vertices, edges and faces of $G$.  By Euler's formula, 
$$
F\ =\ E-V+2-2g\ \geq\ n(k/2-1)+2-2g.  
$$
(the inequality becomes an equality if all faces are simply-connected).  
It follows easily that if infinitely many simple $k$-regular 
graphs can be embedded in a fixed surface $S$, then $k\leq 6$.  
We denote by $\cX_{n,k,g}$ the set of all simple 
$k$-regular graphs embeddable on the orientable surface $S_g$.  

McKay proves in \cite{McKay:83} that for a random $k$-regular graph on 
$n$ vertices, $(\tauG)^{1/n}\to\sigma_k=(k-1)^{k-1}/(k(k-2))^{k/2-1}$ 
as $n\to\infty$.  
We now state the main result of this section: 
\begin{theorem}\label{thm:6sur}
Let $\cX_{n,k,g}$ be as above.  Then 
\begin{equation}\label{sur6:trees}
\lim_{n\to\infty}\left(\sup_{G\in \cX_{n,6,g}} \tauG^{1/n}\right)\ \leq\ 
\sigma_3^2\ <\ \sigma_6.
\end{equation}
\end{theorem}

To prove the theorem we first establish a proposition estimating the 
ratio of the 
tree number of $G$ to the tree number of its dual graph $G^*$. We next  
remark that the dual of a $6$-regular graph on $n$ vertices is ``roughly'' 
a $3$-regular graph on $2n$ vertices; an application of McKay's theorem 
finishes the proof.  

\begin{prop}\label{tree:ratio}
Let $\{G_j\}$ be a sequence of simple graphs on an orientable surface $S$ 
of genus $g$ such that 
their maximal degrees are bounded.  Let $\tauG$ and $\tau(G^*)$ 
be the corresponding numbers of spanning trees, and let $n_j\to\infty$ be 
the number of vertices of $G_j$.  Then 
\begin{equation}\label{tree:surfbound}
n_j^{-2g}\ \ll\ \frac{\tau(G_j)}{\tau(G_j^*)}\ \ll\ n_j^{2g}.
\end{equation}
\end{prop}

\bld{Proof.}
We first define a mapping from the spanning trees 
of $G_j$ into (subsets of) the spanning trees of $G_j^*$ as follows: take 
a spanning tree $T$ of $G_j$, and consider the complement $S=G_j-T$ of $T$ 
in $G_j$.  It is easy to see (\cite{Biggs:dual}) that $S$ is connected, 
spans $G_j^*$ and that one can get a spanning tree of $G_j^*$ by removing at 
most $2g$ edges of $S$.  We denote by $F(T)$ the set of all spanning trees 
of $G_j^*$ which can be obtained from $T$ by the above procedure.  By 
reversing the construction, it is easy to see that the mapping 
$T\to F(T)$ is onto.  

Let $|V(G_j)|=n_j$, and let the maximal degree be bounded above by $\delta$.  
Let $T$ be a spanning tree of $G_j$, and let $S=S(T)=G_j-T$.  Then 
$|S|=|E(G_j)|-n_j+1\leq n_j(\delta/2-1)+1$ and so 
$$
|F(T)|\ \leq\ {|S(T)|\choose 2g}\ \ll\  n_j^{2g}, 
$$
Adding the inequalities for all $T$ proves one of the bounds in 
\eqref{tree:surfbound}.  The other inequality is proved similarly.  
We remark that the constants in \eqref{tree:surfbound} depend on 
$\delta$, $k$ and $g$ only.  

\qed

Consider now a simple $6$-regular graph $G$ (embedded on a surface $S_g$) 
with $n$ 
vertices and $m=3n$ edges, and let $G^*$ be the dual graph 
(of genus $h\leq g$).    
$G^*$ has $m'=3n$ edges and at least $2n+2-2h\geq 2n+2-2g$ vertices; 
moreover, 
since $G$ is simple, every face of $G$ has at least 
three sides, so the minimal degree of $G^*$ is at least three.  It is easy 
to see 
that $G^*$ can be obtained from a (possibly disconnected) graph $G_1$ of 
maximal 
degree at most $3$ with at most $2n+2$ vertices by adding $2h$ ``new'' edges 
between the vertices of $G_1$.  

\bld{Proof of Theorem \ref{thm:6sur}.}
Let $G$ be a graph in $\cX_{n,6,g}$.  By Proposition \ref{tree:ratio}, it 
suffices to estimate $\tau(G^*)$.  Consider a spanning tree $T$ of $G^*$ 
containing $l$ of the ``new'' edges $G^*-G_1$, where 
$0\leq l\leq 2h\leq 2g$.  Removing the $l$ edges 
from $T$ produces a spanning forest of $G_1$ with $l\leq 2h$ components.  
It suffices to establish the estimate \eqref{sur6:trees} for $l$-component 
spanning forests of $G_1$ for each $l$ separately.  For $l=0$ the estimate 
follows easily from McKay's estimate \eqref{limsup:trees} for cubic graphs 
on at most $2n+2$ vertices.  

To establish the estimate for $l\geq 1$, we first remark that the number of 
the connected components of $G_1$ is at most $2h\leq 2g$; let 
$H_1,\ldots ,H_p$ be denote these components, and let $H_k$ have $q_k$ 
vertices, where $q_1+\ldots+q_p\leq 2n+2$. 
Each component $H=H_k$ is a graph of maximal degree at most $3$.   

Let $F$ be an $l$-component spanning forest of $G_1$.  Let 
$F\cap H_k=F_k$ be a spanning forest of $H_k$ with $l_k$ components,   
where $l_1+\ldots+l_p=l$.  $F_j$ can be obtained from a spanning 
tree of $H$ by removing $l_k$ edges from the tree.  Each tree has $q_k-1$ 
edges, so one can get at most ${q_k-1\choose l_k}$ spanning forests by 
removing $l_k$ edges from it, so the total number of such forests is at most 
$$
\tau(H_k)\; {q_k-1\choose l_k}
$$
Multiplying the inequalities for all $H_k$-s we estimate the number of 
$l$-component spanning trees of $G_1$ by 
\begin{equation}\label{tree:surfest1}
\prod_{k=1}^p\tau(H_k)\;\prod_{k=1}^p {q_k-1\choose l_k}
\end{equation}
Now, let a sequence $G_j\in\cX_{6,n_j,g}$ be a sequence of $6$-regular 
graphs on $S_g$, and let $n_j\to\infty$.  Taking the $n_j$-th root 
in \eqref{tree:surfest1} and letting $n_j\to\infty$ shows that it suffices to 
prove that 
\begin{equation}\label{tree:surfest2}
\limsup_{n_j\to\infty}\left(\prod_{k=1}^{p_j}\tau(H_{k,j})\right)^{1/n_j}\ 
\leq\ \sigma_3^2
\end{equation}
where the subscript $j$ refers to the graph $G_j$.  

In order to arrive at a contradiction, we shall assume  
that there exists a sequence of graphs $G_j$ 
on $S_g$ such that the limit in \eqref{tree:surfest2} is strictly 
greater than $\sigma_3^2$.  Passing to a subsequence if necessary we 
may assume that the number $1\leq p_j\leq 2g$ of components of $G_{1,j}$ 
is actually constant and equal to $p$.  We shall order the components 
$H_{k,j}$ so that $q_{k,j}$-s are non-decreasing.  

Passing once again to a subsequence if necessary, we may assume that 
there exists a nonnegative number $r<p$ and a sequence 
$1\leq q_1\leq q_2\leq\ldots\leq q_r$ such that 
$q_{k,j}=q_k$ for $k\leq r$, and $q_{k,j}\to\infty$ for $k>r$.  
Now applying McKay's bound \eqref{limsup:trees} for each $k>r$ 
(with $n=q_{k,j}$) and estimating each of the first $r$ terms in 
\eqref{tree:surfest2} by an absolute constant leads easily to a 
contradiction, finishing the proof of Theorem \ref{thm:6sur}.
  
\qed

The bound proved in Theorem \ref{thm:6sur} is not optimal.  We 
define the optimal bound for $3\leq k\leq 6$:
\begin{equation}\label{treesup:surf}
   \lim_{n\to\infty}\; \left(\sup_{G\in \cX_{n,k,g}} (\tauG)^{1/n}\right)
   \ :=\ \sigma_{k,g} 
\end{equation}
It follows from the proof of Theorem \ref{thm:6sur} that for $g_1<g_2$, 
$$
\sup_{G\in \cX_{n,k,g_2}} \tauG\ \ll\ n^{2g_2} 
\left(\sup_{G\in \cX_{n,k,g_1}} \tauG\right). 
$$
It follows that the limit in \eqref{treesup:surf} does not depend on $g$ 
(as long as $\cX_{n,k,g}$ is not empty for large $n$, which can happen for 
small $g$).  Accordingly, we define $\sigma_s(k)$ to be\footnote{Here 
the subscript $s$ stands for {\em surface}.} the limit 
$\sigma_{k,g}$ achieved in \eqref{treesup:surf} for some fixed $g\gg 0$ and 
$3\leq k\leq 6$.  Noga Alon has informed the authors that he 
can show that 
$$
\sigma_s(k)\ <\ \sigma_k  
$$
(where $\sigma_k$ is defined in \eqref{limsup:trees}).  

It follows from the proof of Theorem \ref{thm:6sur} that 
$\sigma_s(6)\leq\sigma_s(3)^2$.  One can actually show by a similar 
argument that 
$$
\sigma_s(6)\ =\ \sigma_s(3)^2
$$
We make a conjecture about the graphs that achieve the limit in 
\eqref{treesup:surf} for $k=3,4$:
\begin{conj}\label{conj:surftree}
The limit $\sigma_s(3)$ (respectively, $\sigma_s(4)$) is achieved by 
quotients of the hexagonal lattice (respectively, the square lattice) 
in $\reals^2$.  
\end{conj}
  
It would follow from Conjecture \ref{conj:surftree} 
(cf. \cite[pp. 246--247]{CDS}) that 
\begin{equation}\label{treesurfexp:4}
\sigma_s(4)=e^{4G/\pi}\approx 3.21 <  3.375 = \sigma_4
\end{equation}
Here $G$ is Catalan's constant.  The exponent $4G/\pi$ in 
\eqref{treesurfexp:4} is equal to the toplogical entropy 
$$
\int_0^1\int_0^1\log(4-2\cos(2\pi x)-2\cos(2\pi y))\; dxdy
$$
of the {\em essential spanning forest process} on the nearest 
neighbor graph on $\zed^2$ (\cite{Solomyak}).  The integral above 
is equal to the {\em Mahler measure} of the polynomial 
$4-x-1/x-y-1/y$; it is equal to 
$$
2L^{'} (\chi,-1)
$$
where $L(\chi,s)$ is the Dirichlet $L$-function associated to 
the quadratic extension $\ratls(\sqrt{-3})/\ratls$ (cf. \cite{RV}).  
It would be interesting to find candidates for achieving $\sigma_s(5)$.

\subsection{Jacobians of regular graphs}
\label{sec:lattice}
In this section we want to establish an analog for graphs of a result 
of Buser and Sarnak (\cite{Buser:Sarnak}) about the rate of growth of 
the minimal norm of a period matrix of a Riemann surface.  They have 
shown that as the genus of a surface increases, the minimal norm grows 
at most logarithmically in the genus.  

Bacher, de la Harpe and Nagnibeda in their work (\cite{BHN,Nag}) 
defined two lattices associated to any unweighted graph, 
{\em the lattice of integral flows} $\Lambda^1$ and 
{\em the lattice of integral coboundaries} $N^1$.  The quotient 
$(\Lambda^1)^\#/\Lambda^1$ is a finite abelian group called a 
{\em Jacobian} of a graph (in an analogy to Jacobians of Riemann surfaces).  
 
We shall now consider $k$-regular graphs for some fixed $k\geq 3$.  
Let $G\in\cX_{n,k}$; the dimension of $\Lambda^1(G)$ is equal to 
$kn/2-n+1$.  The determinant of $\Lambda^1(G)$ is 
equal to the tree number $\tauG$, and the (unnormalized) minimal norm 
$\tilde{\nu}(G)$ of $\Lambda^1(G)$ is equal to the girth of $G$.  

We rescale $\Lambda^1$ so that the volume of its fundamental domain is equal 
to $1$; the (normalized) minimal norm $\nu(G)$ is given by    
$$
\nu(G)\ =\ \gamma(G)\;\tauG^{2/(kn/2-n+1)}.
$$  
 From \eqref{limsup:trees} and \eqref{liminf:trees} we conclude that 
 for every $\varepsilon>0$,  
$$
\delta_k^{2/(k/2-1)}-\varepsilon\ < \ \tauG^{2/(kn/2-n+1)}\ < \ 
\sigma_k^{2/(k/2-1)}+\varepsilon
$$
for $n$ large enough.  

We want to study the behavior of $\nu(G)$ when $k$ is fixed and 
$n\to\infty$.  Using the known results about the behavior of the girth, 
we conclude: 
\begin{prop}\label{jacobian:rate}
The minimal norm $\nu(G)$ of the Jacobian of a $k$-regular graph 
$G\in\cX_{n,k}$ satisfies 
$$
\nu(G)\ =\ O(\log_{k-1} n)
$$ 
\end{prop}  
This establishes an analogue of the results of  Buser and Sarnak.


\section{Some explicit constructions of cubic graphs}
\label{sec:construct} 

It is known (\cite{BG,Kotzig,Tutte}) that every $3$-connected 
cubic graph may be constructed from a tetrahedron by a sequence of 
{\em edge insertions}: given a cubic graph on $2n$ vertices, choose a pair 
of edges, subdivide them by putting a single vertex in the middle, and 
connect the two subdivision vertices to get a cubic graph on $2n+2$ 
vertices.  So, every $3$-connected cubic graph $G$ may be 
constructed from the tetrahedron provided we know the ``history,'' i.e. the 
sequence in which insertions were applied to successive pairs of edges.  
Clearly, that ``history'' is not unique.  However, one may hope that if $G$ 
is an extremal graph, among all its histories there exist ones with some 
nice properties.  Below we shall give some examples when that is indeed 
the case.   

The first example is the following sequence $S$: starting 
from a tetrahedron, insert edges in such a way that the number of spanning 
trees is maximized after each application.  When the number of vertices is 
small, one may find the graphs in $S$ from the table of graphs in 
\cite{Valdes}.  The $6$-vertex graph is $K_{3,3}$, the $8$-vertex graph is 
a M\"obius ladder (an $8$-cycle with the pairs of ``opposite'' vertices   
connected, cf. \cite[\S 3e]{Biggs:book}); it has $392$ trees.  The $10$-vertex 
graph is the Petersen graph $P$ ($2000$ trees), and the $12$-vertex graph is 
obtained from $P$ by applying $A$ to a pair of edges at distance $2$ from 
each other (this is graph \# 84 in \cite{Valdes}; it has $9800$ spanning 
trees). The last graph also appears in \cite{Holton:cycles}.  From the table 
of graphs in \cite{Valdes} one can check that the graphs in the sequence 
have the most trees among the cubic graphs with the same number of vertices.

The rest of the graphs in the sequence $S$ are shown in the Appendix; they   
were obtained by a computer program. It is remarkable that the sequence  
contains all $(3,g)$-cages for $3\leq g\leq 8$.  Namely, the $14$-vertex 
graph in $S$ is the Heawood graph, the $24$-vertex one is the McGee 
graph, and the $30$-vertex one is the Tutte-Coxeter graph (cf. \cite{Wong} 
or \cite[Ch. 6]{HS:book}).  Unfortunately, the graphs in $S$ with $32$ or 
more vertices seem to be less interesting.  The first few graphs in the 
sequence appear in \cite{AS} (but without the reference to the tree number); 
however, the inclusion of the $12$-vertex graph in the sequence and the 
graphs with $16$ or more vertices appear to be new.  The authors hope that 
similar constructions may lead to discovery of new extremal 
graphs.\footnote{Other sequences of cubic graphs might be constructed  
using other graph invariants (e.g. mean distance -- cf. \cite{CCMS} --  
instead of the tree number) 
in the algorithm; one can also start from a ``large'' cubic graph and delete 
an edge together with its endpoints so that a given invariant is maximized. 
For regular graphs of degree $k\geq 4$ one can use operations similar to 
the edge insertion to increase the number of vertices and construct 
sequences of $k$-regular graphs by maximizing a certain invariant.}  

Edge insertion has its disadvantages: for example, certain 
cages of even girth have the property that every two edges are contained in a 
common shortest cycle (\cite{KT}); inserting an edge will 
then certainly decrease the girth.  So, one might 
try adding several new vertices and edges simultaneously.  For cubic graphs, 
one can try subdividing several edges of an ``old'' graph (with one, two or 
more points on an edge) to get an {\em even} number of ``new'' vertices 
of valence $2$, then choosing a perfect matching for the new vertices and  
connecting pairs of these vertices according to the perfect matching to get 
a new cubic graph.  

Several cubic cages can be obtained by taking uniform subdivisions 
of smaller graphs and inserting edges:\footnote{Similar examples appear 
in \cite{AS}.}
\begin{itemize} 
\item[$\bullet$] The Petersen graph can be constructed from a tetrahedron by 
subdividing its $6$ edges;
\item[$\bullet$] The Heawood graph can be obtained from the unique loopless 
two-vertex cubic graph ($G_2$) by putting 4 new vertices on each of 
its $3$ edges;
\item[$\bullet$] The Tutte graph can be obtained from $2K_{3,3}$ by
subdividing  each of its $18$ edges.

\end{itemize}


\section*{Appendix: a sequence of cubic graphs; constructing cubic cages.}
\label{sec:sequence} 

This section describes a sequence of cubic graphs generated (by a computer) 
from the tetrahedron by edge insertion where at each stage the insertion was 
chosen so that the number of spanning trees of he new graph was the maximal 
possible.  The first six graphs are shown in Figure 1; the endpoints of the 
new edge inserted into a cubic graph on $2n$ vertices were numbered $2n+1$ 
and $2n+2$.   

\vskip 10pt
\hbox to \hsize
{
\hfill
\epsfysize=3in
\epsffile{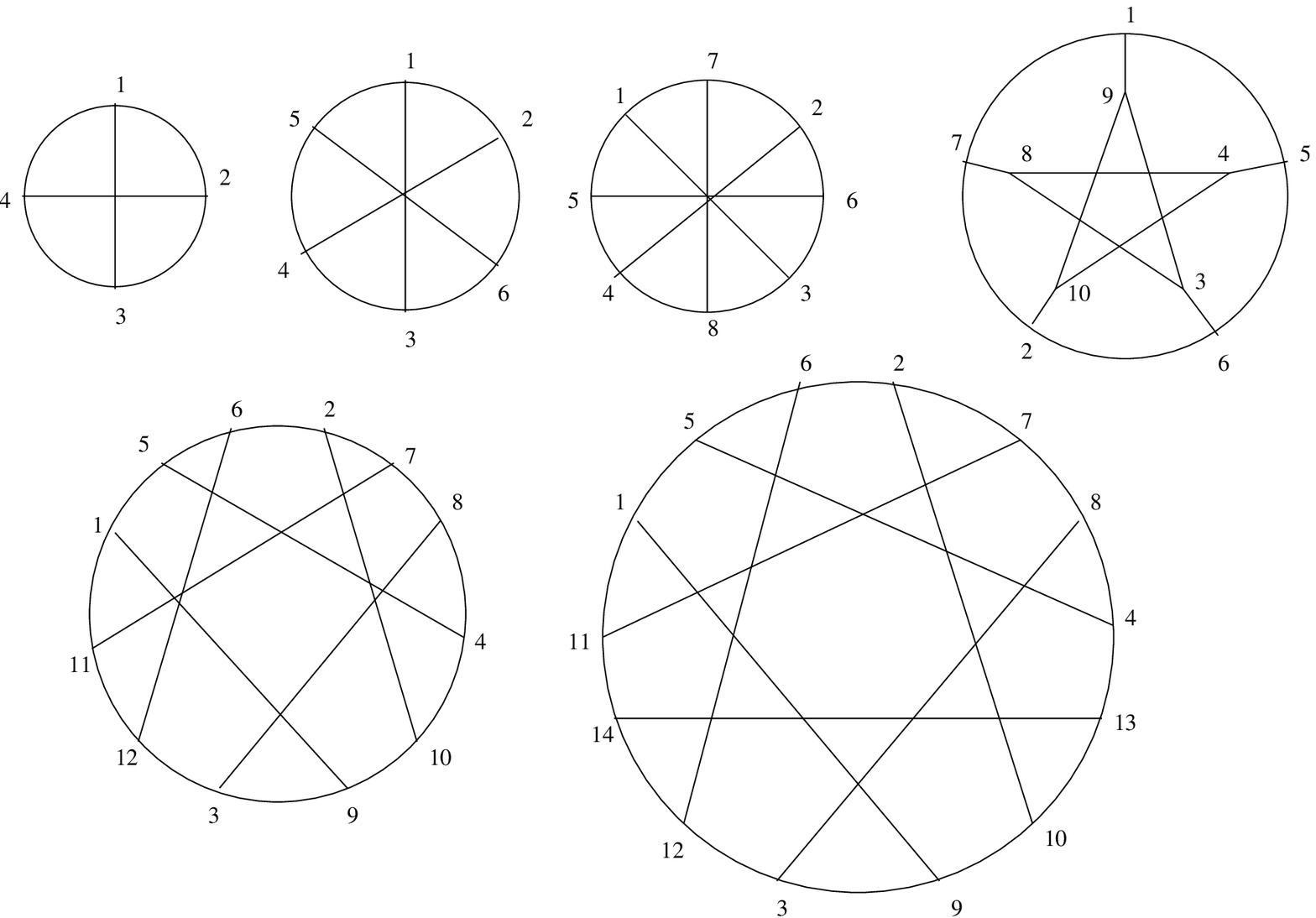}
\hfill
}
\hbox to \hsize
{
\hfill
\hbox to 0pt{\hss Figure 1. The first six graphs in a sequence $S$.\hss}
\hfill
}
\vskip 10pt

The graphs with $18$ and $24$ vertices (as they appear in $S$) are 
shown in Figure 2.  

\vskip 10pt
\hbox to \hsize
{
\hfill
\epsfysize=2in
\epsffile{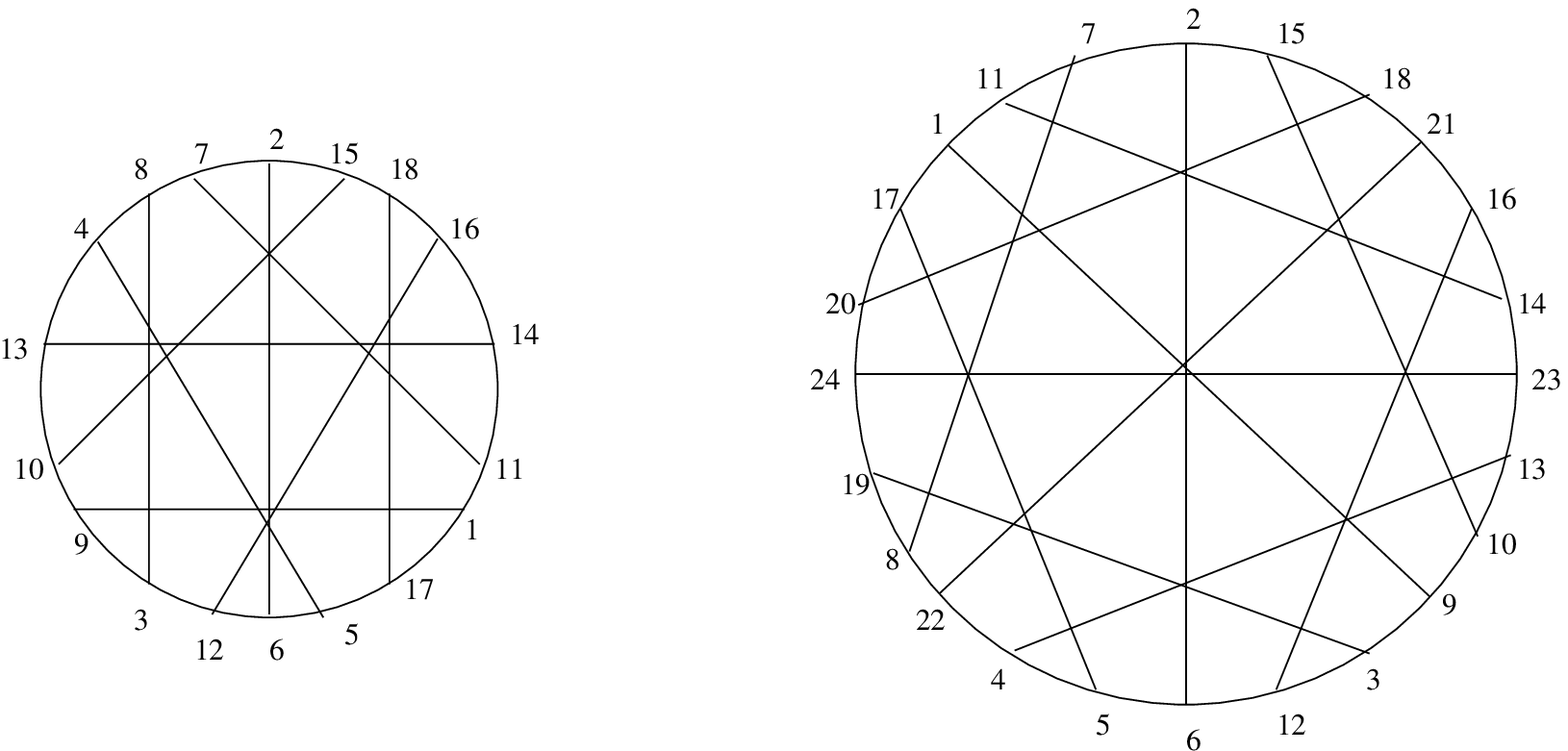}
\hfill
}
\hbox to \hsize
{
\hfill
\hbox to 0pt{\hss Figure 2. The $18$-vertex graph and the McGee graph.\hss}
\hfill
}
\vskip 10pt

The Tutte-Coxeter graph (as it appears in $S$) is shown on Figure 3.

\vskip 10pt
\hbox to \hsize
{
\hfill
\epsfysize=3in
\epsffile{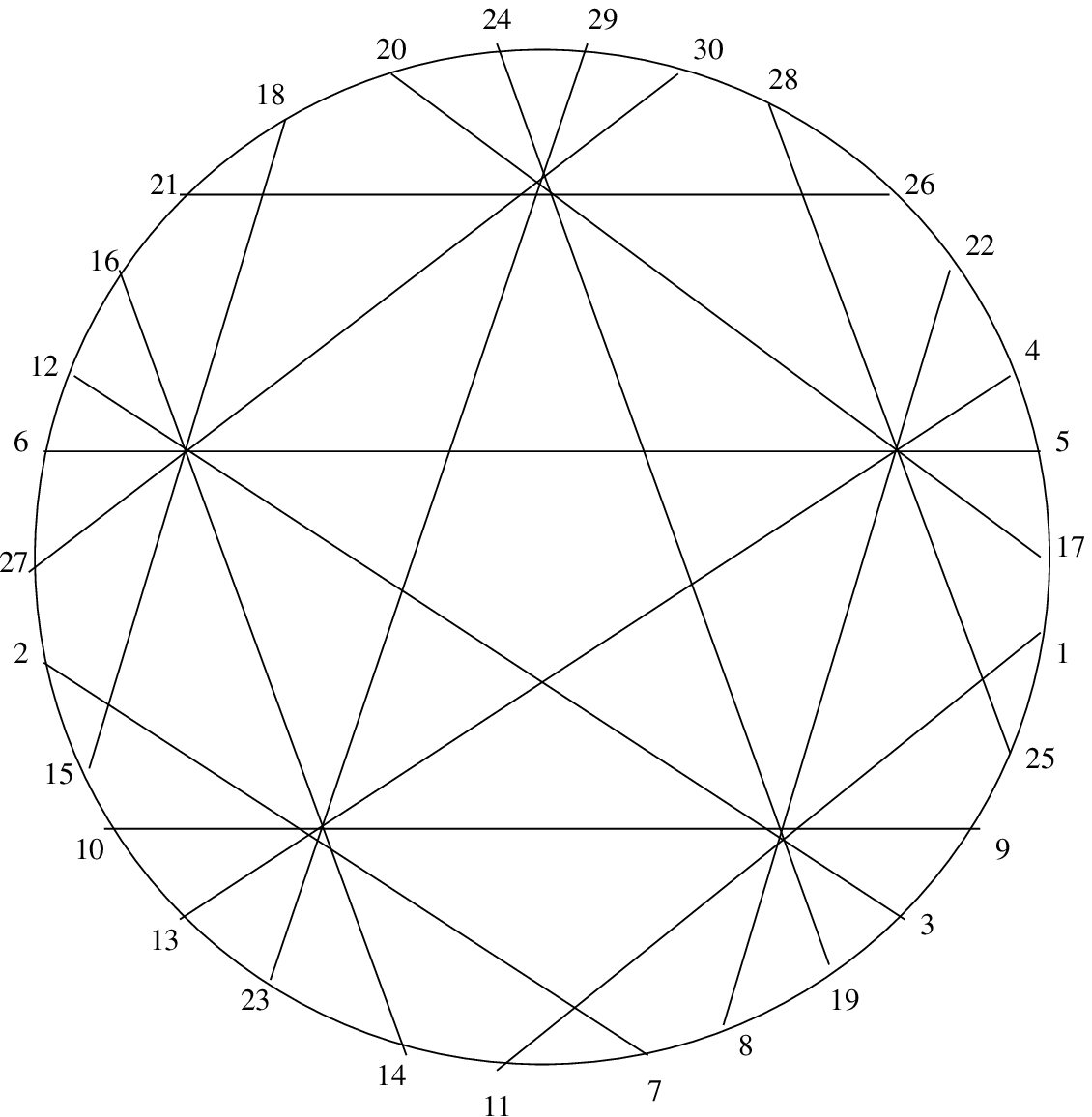}
\hfill
}
\hbox to \hsize
{
\hfill
\hbox to 0pt{\hss Figure 3. The $(3,8)$-cage.\hss}
\hfill
}
\vskip 10pt

In Figure 4 we show how to construct Petersen and Heawood graphs by 
edge insertions into the tetrahedron and the graph $G_2$ respectively.  
The vertices of $K_4$ and $G_2$ are shown in black.

\vskip 10pt
\hbox to \hsize
{
\hfill
\epsfysize=2in
\epsffile{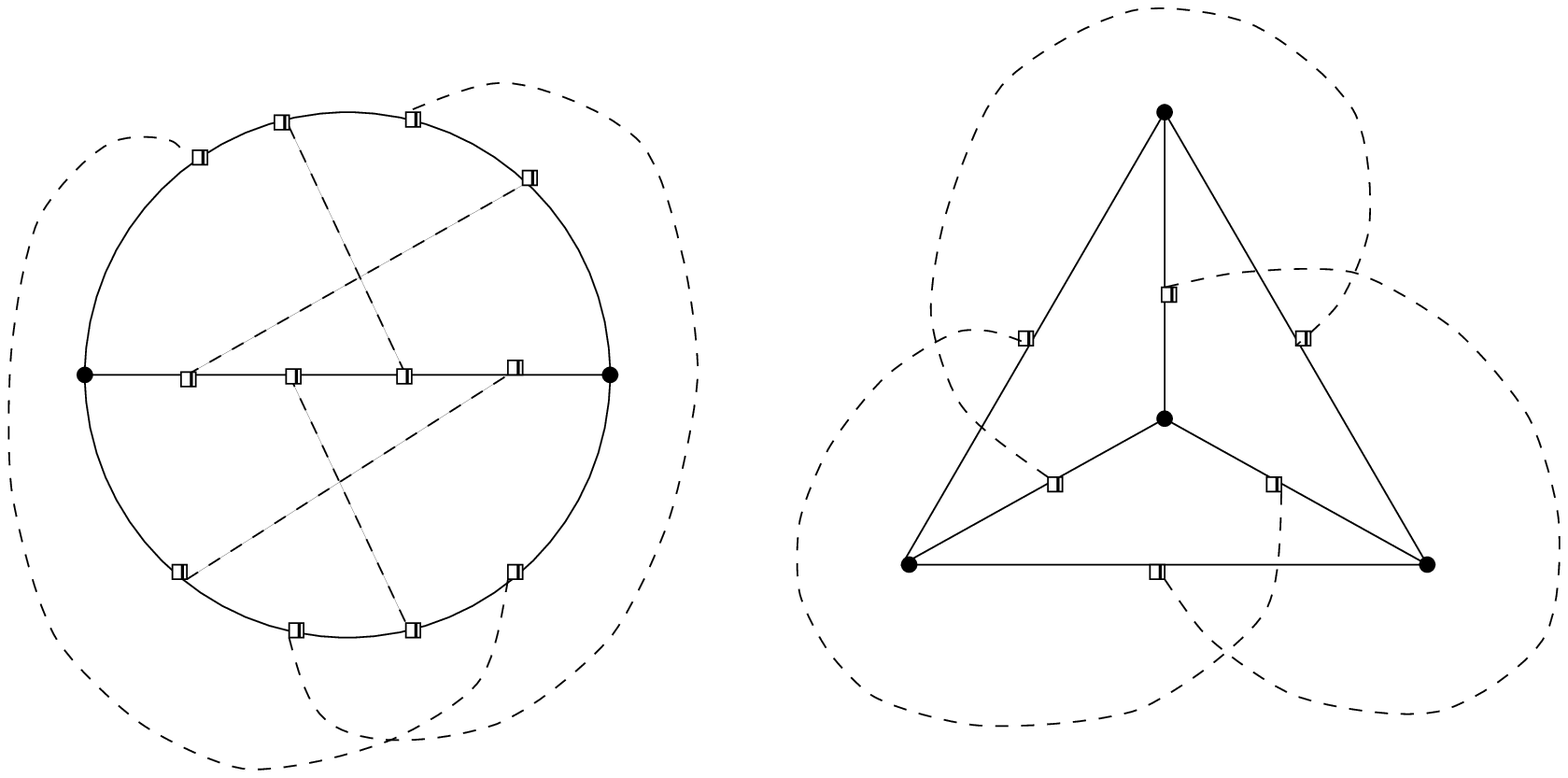}
\hfill
}
\hbox to \hsize
{
\hfill
\hbox to 0pt{\hss Figure 4. Petersen and Heawood graphs.\hss}
\hfill
}
\vskip 10pt

In Figure 5 we show how to construct the Tutte graph by edge insertions into 
$2K_{3,3}$.  The vertices and the edges of $2K_{3,3}$ are shown in black; 
the vertices of the two copies of $K_{3,3}$ are marked $A$ and $B$ 
respectively.

\vskip 10pt
\hbox to \hsize
{
\hfill
\epsfysize=3in
\epsffile{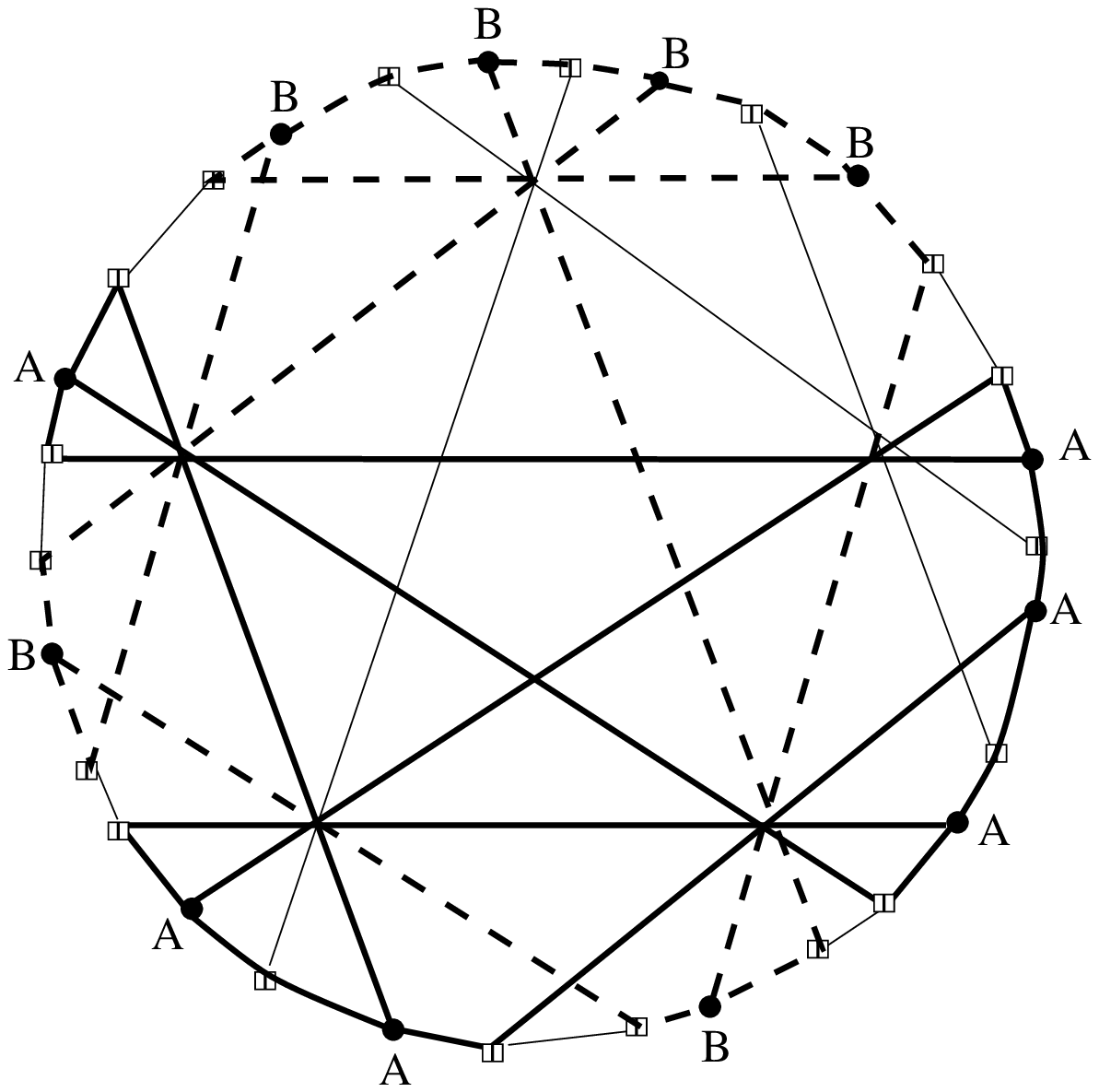}
\hfill
}
\hbox to \hsize
{
\hfill
\hbox to 0pt{\hss Figure 5. The Tutte graph obtained from $2K_{3,3}$.\hss}
\hfill
}

\end{document}